\documentclass[review]{elsarticle}

\usepackage{amsmath,amsfonts,amssymb,fancybox,bm}
\usepackage[T1]{fontenc}
\usepackage{colortbl}
\usepackage{xcolor}
\usepackage{pifont}
\usepackage[nointegrals]{wasysym}
\usepackage{slashbox}
\usepackage{multirow}
\usepackage{hyperref}
\usepackage{algorithmicx}
\usepackage[section]{algorithm}
\usepackage{algpseudocode}
 
\usepackage{makecell}
\usepackage{float}
\usepackage{nomencl}
\makenomenclature
\usepackage{comment}
\usepackage{lineno}
\modulolinenumbers[1]

\DeclareMathOperator*{\argmax}{arg\,max}
\DeclareMathOperator*{\argmin}{arg\,min}
\DeclareMathOperator*{\datasample}{data\,sample}

\journal{Applied Mathematical Modelling}

\begin{document}

\begin{frontmatter}

\title{An Integrated Optimization Model for the Multi-Port Stowage Planning and the Container Relocation Problems}


\author[aut1]{Catarina Junqueira\corref{mycorrespondingauthor}}
\cortext[mycorrespondingauthor]{Corresponding author}
\ead{catarinajunqueira2@gmail.com}
\author[aut2]{Miguel Paredes Qui\~nones}
\author[aut3]{Anibal Tavares de Azevedo}
\author[aut3]{Cleber Dami\~ao Rocco}
\author[aut1]{Takaaki Ohishi}

\address[aut1]{University of Campinas, School of Electrical and Computer Engineering, Campinas, Brazil.}
\address[aut2]{IBM Research, S\~ao Paulo, Brazil.} 
\address[aut3]{University of Campinas, School of Applied Sciences, Limeira, Brazil.}

\begin{abstract}
This paper introduces a new optimization model that integrates the multi-port stowage planning problem with the container relocation problem. This problem is formulated as a binary mathematical programming model that must find the containers' move sequence so that the number of relocations during the whole journey of a ship, as well as the associated port yards is minimized. Modeling by binary variables to represent the cargo status in a ship and yards makes the problem very complex to be solved by exact methods. To the best of our knowledge, this integrated model has not been developed yet as that such problems are always addressed in a partitioned or hierarchical way. A demonstration of the benefits of an integrated approach is given. The model is solved in two different commercial solvers and the results for randomly generated instances are presented and compared to the hierarchical approach. Two heuristics approaches are proposed to quickly generate feasible solutions for warm-starting the model. Extensive computational tests are performed and the results indicate that the solution approaches can reach optimal solutions for small sized instances and good quality solutions on real-scale based instances within reasonable computation time. This is a promising model to support decisions in these problems in an integrated way.
\end{abstract}


\end{frontmatter}


\section{Introduction} \label{sec:introduction}

Maritime transport is the most important mode of transport for international trade as more than four fifths of world's merchandise trade by volume is transported by sea \citep{Unctad2019}. Its continuous growth, as well as the up-sizing of container ships, represents greater pressure on port terminals. To maintain competitiveness and productivity, as the organization of the container terminals significantly affects the entire supply chain, it is crucial for ports to adapt to the growing complex issues.

Taking this into account, this paper proposes a binary linear optimization model that integrates the Multi-Port Stowage Planning Problem (MPSP) in the ship with the Container Relocation Problem (CRP) at the port yard, called MPSP-CRP here. According to \cite{Steenken2004}, these two problems can be optimized adopting operations research methods, and are defined as follows:


\begin{itemize}
\item \textbf{Multi-Port Stowage Planning Problem}: This consists of determining the position of the containers on board a ship along its route with the objective of minimizing the number of relocations in the loading and unloading operations at each port \citep{Torres2019}.

\item  \textbf{Container Relocation Problem}: This problem deals with a given set of homogeneous containers stored in a set of two-dimensional last-in-first-out (LIFO) stacks, whose relocations are necessary to retrieve the containers from the stacks while minimizing the number of those relocation \cite{Caserta2012}.
\end{itemize}

A relocation occurs while retrieving a container that is not at the top of a container stack, and is considered to be an unproductive move. The number of relocations is an important factor that affects the operational efficiency in container terminals, and it is correlated with the loading sequence and the re-handling strategy \cite{Ji2015}. 

There are studies in the literature that highlight the importance of reconciling conflicting interests of different port stakeholders to improve overall performance in the whole chain, such as \cite{Ha2019} and \cite{Brooks2013}. 


Specifically in this study, from the perspective of the shipping line, to define the stowage planning knowing the location of the containers at the port yards means ensuring fast and precise loading and unloading operations. Meanwhile, from the perspective of the container terminal managers, integrating seaside and yard decision problems at container terminals make sense in order to plan their resources for transferring cargo both at the seaside and landside operations.

Although port integration has been a topic of research for years, the practice of doing separate or simplified planning still tends to exist. Some examples are: focusing on the export flow at the one-port ship loading problem and ignoring the cargo arrangement in stacks at the yard \cite{Monaco2014, Iris2015, Iris2018}; considering only the cargo arrangement on the ship \cite{Torres2019}; observing only the arrangement in the yard \cite{Rahman2016}; or considering the transshipment flow but ignoring the stacking structure in the ship and the yard \cite{ZHEN2016700}.




Regarding the aforementioned, the MPSP-CRP defines the sequence of container removal from all the port yards belonging to the ship's route so that the number of relocations are minimized both in the yards and on the ship, which is pointed out as the key contribution in this paper.




By integrating these two problems, the impact of a decision in one port can be estimated at the following ones, which will be advantageous for both the port terminals and the shipping lines, so that one will not gain over the detriment of the other. Thus, it can be guaranteed that the stowage planning obtained is executable by each port terminal and that the sequence of container removals from the yards is also good for the container ship.




To the best of our knowledge, there is no mathematical model in the literature that represents these activities in an integrated way, or a methodology that solves it analytically aiming for optimization. Therefore, this paper aims to demonstrate that the CRP at the port yard and the MPSP in the ship are actually correlated and should be optimized jointly.

Here, also a heuristic approach is also presented to quickly generate solutions for the MPSP-CRP and they are analysed with the results obtained from two commercial solvers: Cplex\texttrademark{} and Gurobi\texttrademark{} as exact methods (Branch and Cut). The heuristic solutions are used for warm-starting the exact method in the MPSP-CRP model. Extensive computational tests are performed to demonstrate the complexity in solving problems modeled by the commercial optimization solvers. 




This paper is organized as follows, besides this introduction. Section \ref{sec:description} describes each problem (MPSP and MPSP) and presents the literature review. The integration of the problems is discussed using a mathematical demonstration concerning the benefits of the integrated approach with respect to solving each problem hierarchically. Section \ref{sec:MPSP-CRP-model} presents the proposed optimization model, as well as its assumptions and mathematical formulation. In Section \ref{sec:solution-approach}, the heuristics approaches with the exact method are presented. In Section \ref{sec:computational-results}, the computational results are presented and discussed. Finally, in Section \ref{sec:conclusions} the findings are summarized and some extensions of this research are outlined.

\section{Description of Problems} \label{sec:description}

This section describes the MPSP and CRP problems. To make the explanation clearer, each problem will be described separately, and then the integration between them is justified.

\subsection{The Multi-Port Stowage Planning Problem}

\begin{figure*}[!h]
  \centering
  \includegraphics[width=26pc]{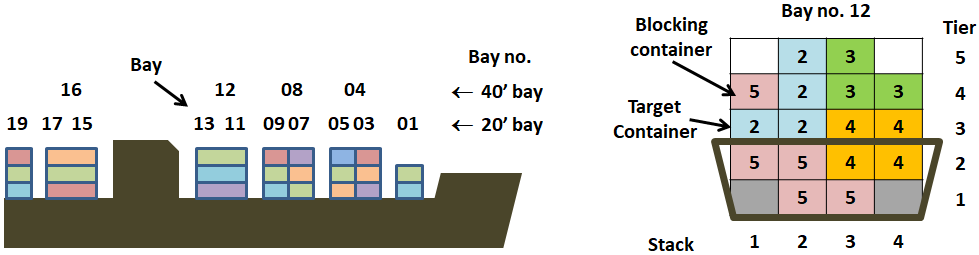}
  \caption{Cellular structure of a container ship.  Adapted from \citep{Bierwirth2010}. \label{figg1}}
\end{figure*}

The formulation of the stowage planning problem is related to the cellular structure that the container ship has, as it can be seen in Figure \ref{figg1}. This structure entails that a container may only be moved if there are no other containers above it. Otherwise, these blocking containers shall be removed to allow the access of the target container. This movement is known as relocation, and it must be minimized to improve port efficiency. Relocation movements can occur frequently and lead to a longer turnaround time of the ship in the berth. To avoid such inconveniences, the stowage plan is prepared so that the decision in one port does not entail many relocations in the next ports to be visited.

The MPSP problem was proven by \cite{Avriel2000} to belong to the complexity category of NP-Complete, by demonstrating that it can be transformed into the problem of coloring circle graphs. Due to its computational intractability for large instances, heuristic and meta-heuristic methods are predominantly used to obtain good solutions in polynomial time, as can be seen in these approaches in the literature: \cite{Botter1991, Avriel1998, Wilson1999, Wilson2000, Dubrovsky2002, Ambrosino2006, Imai2006, Sciomachen2007, Ambrosino2010, Fan2010}.

Nevertheless, some exact models were presented in \cite{CHEN199568}, \cite{Avriel1998}, \cite{Kang2002}, \cite{Pacino2011}, \cite{Monaco2014} and \cite{Torres2019}. However, the mathematical formulations were limited to small-sized instances. For \cite{Lehnfeld2014}, the main methods to solve the stowage planning problem are Dynamic Programming, heuristic approaches, Integer Programming, Genetic Algorithm, Branch and Bound and Tabu-search. \cite{Wilson2001} also adds the simulation based on probability methods and rule-based expert systems to the list. A survey of the state-of-the-art methods used on the ship loading problem, which is part of the stowage planning problem, is given by \cite{Iris2015}.

A variant on the MPSP problem is the single-port stowage planning problem (SPSP), which according to \cite{Torres2019}, consists of determining the arrangement of containers in the ship at a given port without considering loading containers in subsequent ports.


Some papers have taken into account features of the containers when developing the stowage planning, such as the dimension (standard, 45-footer, high-cube, oversized), the weight class (light, medium, heavy), the type (reefer, open-top), the load (dangerous, perishable), and the port of destination (POD), such as \cite{Monaco2014} and \cite{Iris2018}. In \cite{Meisel2010}, containers are allowed to be internally reshuffled. Stability constraints of the vessel and a hatch cover are not considered, and the containers are distinguished only by their destination.


\subsection{The Container Relocation Problem}

\begin{figure*}[h!]
    \centering
    \includegraphics[width=10pc]{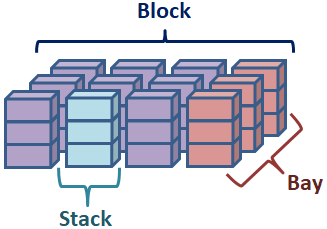}
    \caption{Part of the port yard: block, bays and stacks. Adapted from \citep{JIN2015}. \label{figg2}}
\end{figure*}

A port yard is a place where containers are temporarily stored until they are loaded onto a ship, truck or train (Figure \ref{figg2}). The container relocation problem is considered to be a complex problem because of the uncertainty regarding which container will be needed first and by the limited space of the storage area \cite{Sauri2011}.

The problem of retrieving a set of containers from a given stack in a given order was proven to be NP-Complete by \cite{Caserta2012}. As occurred with the MPSP, this has led researchers to use heuristics and meta-heuristics procedures, as can be seen in \cite{Kim2006, Lee2007, Wan2009, Tang2012, Petering2013}. Exact formulations have been presented in \cite{Wan2009, Tang2012, Petering2013, ExpositoIzquierdo2015, Zehendner2015}; but all the models require a long computational time to solve large instances, since they generate an exorbitant number of variables.

The literature in this problem contains the following main variants: restricted or unrestricted container relocation problems. In the restricted case, the assumption made is that the containers that can be relocated are only those that are blocking the next container to be retrieved. In the unrestricted case, all feasible relocation movements are considered. The consideration of the restricted variant is reasonable as it allows a reduction in the search space, however it can lead to the exclusion of optimal solutions. In our study, the restricted container relocation assumption is adopted.

Both restricted and unrestricted problems are examined by \cite{Zhu2012}, by using iterative deepening A$^{\ast}$ algorithms (IDA$^{\ast}$). The effects of some greedy heuristics are also investigated.



Some studies have also considered the extensions of relocation problems. For example, \cite{Lee2010} and \cite{Silva2018} both considered a bi-objective container relocation problem. \cite{Lee2010} proposed to minimize the number of movements and the amount of time required for crane processing, while \cite{Silva2018} proposed to minimize the number of relocations and the expected number of relocations for retrieving the containers from the next customer.

Another associated issue is the \textit{Pre-Marshalling} problem. It occurs when containers within a storage area must be rearranged so that they can be removed later without any additional relocations. One key feature in the \textit{Pre-Marshalling} is that no containers are removed from the yard and its advantage is that the containers can be removed quickly when needed, since they are rearranged, especially in a short period of time. However, according to \cite{Monaco2014}, in some cases it may be that \textit{Pre-Marshalling} is neither feasible nor profitable. It depends on the time and resources available to reallocate containers before the ship loading begins.

Regarding the \textit{Pre-Marshalling} problem, \cite{Lee2007} present an integer programming model which has a multi-commodity network flow problem embedded within. The authors discuss several possible variations of the model that can be obtained by modifying some of the constraints, and propose a heuristic to solve it more efficiently. \cite{MELODASILVA2018} presented a unified integer programming model for solving the \textit{Pre-Marshalling} and the container relocation problem in the restricted and unrestricted variant. A simple greedy heuristic was implemented for the \textit{Pre-Marshalling} and a Branch and Bound algorithm for the restricted container relocation problem was implemented and used to provide upper bounds for the unrestricted problem.

An overview of the literature on the port yard loading, unloading and \textit{Pre-Marshalling} problems is provided in \cite{Lehnfeld2014}. Some other papers that make surveys on container terminals related problems are \cite{Stahlbock2008}, \cite{Carlo2014}, \cite{Bierwirth2015}, \cite{KimeLee2015} and \cite{LeeeSong2017}. 


The next section explains the integration between the container relocation problem with the stowage planning problem, which is the aim of this paper.

\subsection{Integration}

The idea of integrating decision problems is highly interesting in the literature because it takes into consideration the impact of decisions for: (i) the entire chain of operations through time and not only the local impact for one port, and; (ii) the cooperation between the two most representative agents in port logistics: the shipping line and the terminal operator.

The key feature of these two problems that allow us to work with them together is that both the container ship space and the port yard have a cellular structure and can be represented as matrices in the mathematical models. This structure means that in both problems, cargo should be organized in vertical stacks, and therefore similar characteristics between these problems are derived. One example is the constraint that a container can only be accessed from the top. For that reason, an approach can be developed to jointly solve the stowage planning problem and the container relocation problem. 
The challenge is to choose an exact position for each container among the numerous possible positions in a way that the total number of movements along the entire ship journey is minimized.

Moreover, it could be said that the CRP and the MPSP are indeed strictly correlated as the total number of movements to remove all the containers from the yard and loading them onto their destination ship is just an estimation of the necessary time that the ship will spend docked at the berth. Consequently, the practical implication of integrating these two problems is to increase the profits of both the ship and the terminals as the stowage plan executed in conjunction with the port can lead to a more accurate turnaround time estimate, which could prevent cargo from being left at the port (because there was not enough time) or the ship being fined if it stays longer than it should have moored in the port.  On the port side, fewer movements performed to load the ship means savings in time and money.


For these reasons, it is no coincidence that recent literature indicates that more researchers are attempting to integrate problems that appear in ports, for example: integrated berth allocation and quay crane scheduling \cite{Bierwirth2010, Chang2010, Imai2008, Meisel2009, Yang2012}; integrated allocation of berths and yard operation planning \cite{Hendriks2013}; and allocation of empty containers in the yard coupled with vehicle routing \cite{braekers13}.

In \cite{Izquierdo2013}, some general guidelines for designing integration approaches for the main container terminal decision problems are provided. The authors concluded that "in-depth research of integration proposals within the seaside and yard area logistical problems is essential to increase and maintain a high productivity as well as produce real-context planning solutions."

This recent tendency is justified by the fact that the optimization of operations for just one stage does not increase overall port efficiency, because further and non-optimized stages behave as bottlenecks. Thus, it makes sense to use an integrated approach and solve the CRP and the MSPS, to gain efficiency and to make a more effective use of assets.


To show the benefits of the integrated model with respect to solving each problem hierarchically, let us define the optimal solution of the integrated problem as:

\begin{align}
P_{1}(x^*,y^*) = ~ & \underset{x,y}{\text{Min}}
& c^Tx + d^Ty \label{eq:optProb} \\
& \text{subject to:} & Ax \leq b \label{eq:constraint2} \\
& & By \leq f \label{eq:constraint3} \\
& & Mx + Hy \leq g \label{eq:constraint4}
\end{align}

Where constraints \eqref{eq:constraint2} represent those related to the CRP, constraints \eqref{eq:constraint3} represent those related to the MPSP and constraints \eqref{eq:constraint4} represent the integration constraints.

In a hierarchical solution scheme, the MPSP, represented in problem $P_2$, must be solved first. By doing so, the optimal solution $y^{+}$ is obtained, which represents the stowage planning. This solution is then passed to problem $P_3$ to obtain a solution for CRP ($x^{+}$), which will be simultaneously feasible for the MPSP problem as meeting constraint \eqref{eq:constraint8} will enforce it. As stated in \cite{Caserta2012}, the accessibility of any container at the yard is guaranteed, then any feasible solution $y^+$ into $P_3$ will result in a feasible problem.

\begin{align}
P_{2}(y^+) = ~ & \underset{y}{\text{Min}}
& d^Ty \label{eq:optProb2} \\
& \text{subject to:} & By \leq f \label{eq:constraint6} 
\end{align}

\begin{align}
P_{3}(x^+) = ~ & \underset{x}{\text{Min}}
& c^Tx  \label{eq:optProb3} \\
& \text{subject to:} & Ax \leq b \label{eq:constraint7} \\
& & Mx + Hy^+ \leq g \label{eq:constraint8}
\end{align}

Observe then that the solutions $(x^+,y^+)$ are feasible for $P_1$ since they meet constraints \eqref{eq:constraint2},  \eqref{eq:constraint3} and  \eqref{eq:constraint4}. Therefore, by basic optimization theory \citep{Nocedal2006}, one can conclude that $P_{1}(x^*,y^*) \leq P_{1}(x^{+},y^{+}) \longrightarrow P_{1}(x^*,y^*) \leq P_{2}(y^+) + P_{3}(x^+)$. Meaning that any feasible solution of the hierarchical problem is going to be worse, or in the best case, equal to the optimal solution of the integrated problem (represented by $P_{1}(x^*,y^*)$). In the next sections, we show some insights into the characteristics of the solution of the CRP and MPSP problems. 

Taking this into account, the next section presents the integrated optimization model for the multi-port stowage planning problem and the container relocation problem.

\section{MPSP-CRP Model} \label{sec:MPSP-CRP-model}

\begin{figure}[h!]
\centering
\includegraphics[width=17pc]{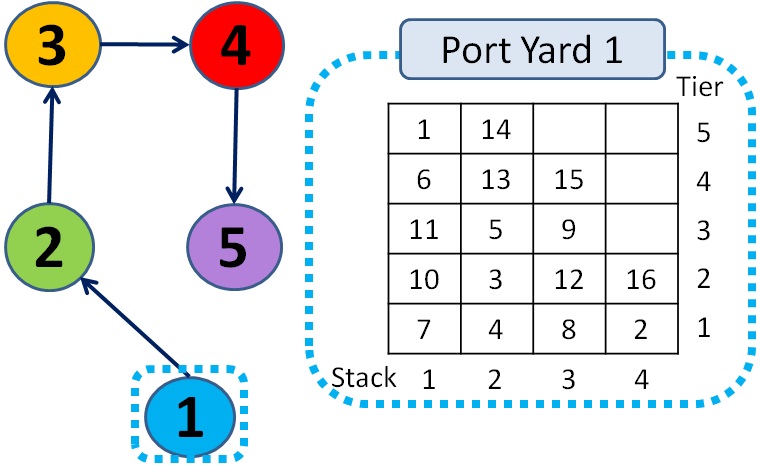}
\caption{Example of containers in Port Yard 1 for a ship route over five ports.}
\label{Descricao_Problema_Rota}
\end{figure}

\begin{figure*}[h!]
\centering
\includegraphics[width=26pc]{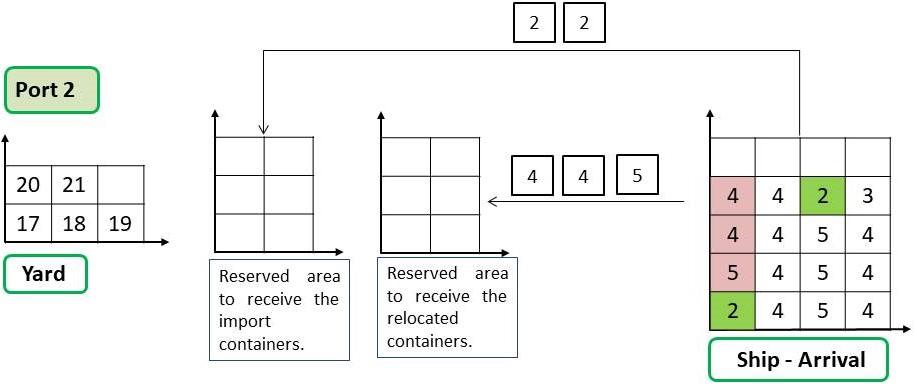}
\caption{Example of ship unloading in the Port Year 2.}
\label{Descricao_Problema_Descarregamento}
\end{figure*}

\begin{figure*}[h!]
\centering
\includegraphics[width=26pc]{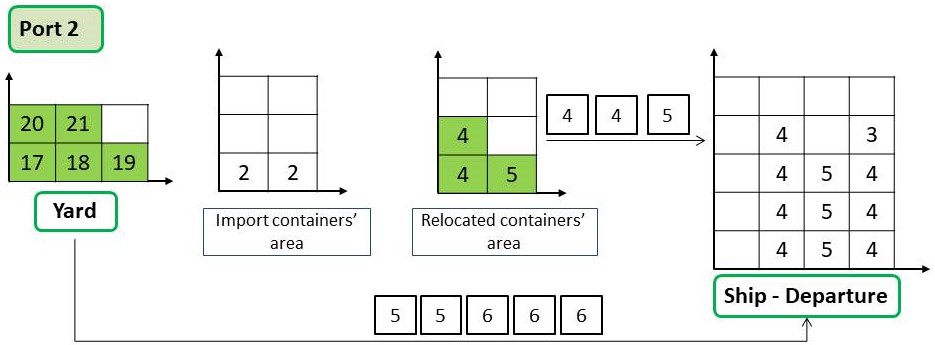}
\caption{Example of ship loading in the Port Yard 2.}
\label{Descricao_Problema_Carregamento}
\end{figure*}

In the MPSP-CRP,a ship's journey is represented by a set of ports $P$. In each port $p \in \{1, ..., P-1\}$ in which the ship docks, there is a yard with containers that must be loaded on it (Figure \ref{Descricao_Problema_Rota}). At the first port $(p = 1)$, the ship arrives empty and is loaded with containers that were in the yard. In the next ports $p \in \{2,...,P-1\}$, the ship first unloads the containers destined to port $p$ where it is, as illustrated in Figure \ref{Descricao_Problema_Descarregamento}, and then receives the loading of containers destined to the following ports (Figure \ref{Descricao_Problema_Carregamento}). Finally, when the ship arrives at the last port of its journey $P$, only the unloading is performed. 

Note that, to remove all containers from each yard belonging to each port $p$, a problem similar to the block relocation problem (BRP) discussed by \cite{Caserta2012} must be solved. The difference is that here no pre-defined removal order is imposed. The only pre-information required by the MPSP-CRP is the location of each container at the yard and their destination. Thus, the MPSP-CRP decides the removal sequence of each container in each port yard by looking at the entire journey of the ship.

It can also be observed that, in Figure \ref{Descricao_Problema_Rota} the containers in the port yard are represented by different numbers in the ship (see Figures \ref{Descricao_Problema_Descarregamento} and \ref{Descricao_Problema_Carregamento}). This is because in the port yard, the number represents only the container, not the removal order or port destination. Meanwhile, in the ship, the numbers represent the destination port of each container.

In Figure \ref{Descricao_Problema_Descarregamento}, it can be concluded that the ship arrives at the second port of its route and it has to unload the two highlighted containers (number 2), of which one of them is not at the top of its stack. For this reason, first the blocking containers (two containers number 4 and one number 5) located above the one target container will be unloaded to a reserved area of the yard. Next, it can see in Figure \ref{Descricao_Problema_Carregamento} that after all containers destined to this port were unloaded from the ship, the container loading operation at the second port yard and the relocated containers is started, replacing the two containers headed for Port 4 and the one headed for Port 5. Moreover, two more containers (number 20 and 21) are also destined to Port 5 and the three containers (number 17, 18 and 19) are destined to Port 6. The MPSP-CRP model follows this representation and it is described in details in the following section.


\subsection{Model assumptions}

The MPSP-CRP extends and integrates the models presented in \cite{Caserta2012} and \cite{Avriel1998}. The following assumptions have been made for the sake of simplicity, without compromising the general application of the solutions: 

\begin{enumerate}[i]

    \item The container ship has a rectangular format and can be represented as a matrix with tiers ($r \in R$) and stacks ($c \in  C$) with maximum capacity of $R \times C$ containers. Each space to store a container is represented by a coordinate ($r,c$), where $r \in \{1, ..., R\}$ and $c \in \{1, ..., C\}$. An irregular format may be achieved by simply adding constraints which represent imaginary containers that occupy the same spaces throughout the whole route \cite{Ding2015}.
    
    \item Each yard $o$ can also be represented as a matrix with tiers ($h \in H_{o}$) and stacks ($w \in W_{o}$). Each space to store a container is a coordinate ($i,j$), where $i\in \{1,... ,W_{o}\}$, $j\in \{1,... ,H_{o}\}$ and $o \in \{1,...,P-1 \}$.
    
    \item All containers have the same size and weight.
    
    \item All containers are accessible only from the top of the stack.
    
    \item The container ship can always carry all the containers in each port and its capacity will never be exceeded.
    
    \item To guarantee accessibility of all containers in each yard $o$, there must be at least $H_{o} - 1$ empty slots, where, $W_{o} \times  H_{o} \geq N_{o} + (H_{o} - 1)$ should be satisfied, where $o \in \{1,...,P-1\}$.
    
    \item Each time period $t \in T_{o}$ is defined by a complete set of moves to retrieve a container from the yard and load it onto the ship, where $o \in P$. This includes relocation moves, if needed.
    
    \item All containers in the port yard will be loaded into the ship.
    
    \item The unloaded containers are destined to a specific area of the port yard and are not mixed with the containers that are going to be loaded, as shown in Figures \ref{Descricao_Problema_Descarregamento} and \ref{Descricao_Problema_Carregamento}.

\end{enumerate}

Note that it is possible to achieve a three-dimension representation by just replicating the proposed representation for the desired number of bays. Additionally, the assumptions made about the containers' weight and size are not a limitation of the model since it enables a focused analysis on how ship and yard container arrangements are affected due to integrated planning.

\subsection{Mathematical formulation}

\subsection*{Sets}

\noindent $P$: ports.\\
\noindent $N$: containers.\\
\noindent $O$: yards.\\
\noindent $R$: tiers in the ship.\\
\noindent $C$: stacks in the ship.\\
\noindent $\varOmega_{o}$: subset of containers in yard $o$.\\
\noindent $W_{o}$: subset of stacks in yard $o$.\\
\noindent $H_{o}$: subset of tiers in yard $o$.\\
\noindent $T_{o}$: time periods in yard $o$ - equal the number of containers in yard $o$.\\
\noindent $\phi_{od}^{n}$: containers $n$ that should be transported from origin $o$ to destination $d$.

\subsection*{Indexes}

\noindent $p,e,m$: indexes for ports, $p,e,m \in \{1,...,P\}$.\\
\noindent $n$: index for containers, $n \in \varOmega_{o}$.\\
\noindent $i,i'$: index for yard's tier, $i,i' \in W_{o}$.\\
\noindent $j,j'$: index for yard's stack, $j,j' \in H_{o}$.\\
\noindent $k$: index for yard's tiers when a relocation is performed, $k \in W_{o}$.\\
\noindent $l,l'$: index for yard's stacks when a relocation is performed, $l \in H_{o}$.\\
\noindent $r$: index for ship's tiers, $r \in R$.\\
\noindent $c$: index for ship's stacks, $c \in C$.\\
\noindent $t,t'$: index for time period, $t \in \{1,...,T_{o}\}$.\\
\noindent $o, o'$: index for port origin for a given container, $o,o' \in \{1,..., P-1\}$.\\ 
\noindent $d$: index for port destination for a given container, $d \in \{o+1,...,P\}$.\\
\noindent $a$: index for the port where the container is relocated, $a \in \{o+1,...,d\}$.

\subsection*{Parameters}

\noindent $N_{o}$: total number of containers in yard $o$.\\
\noindent $\theta_{o}$: total number of containers in the ship when leaving port $o$.\\
\noindent $F_{od}$: number of containers with origin $o$ to destination port $d$, with dimension $(P-1)\times(P-1)$.\\
\noindent $M$: a sufficiently large value (Big-M).

\subsection*{Decision variables}

The variables are divided into three groups. The first group comprises the variables exclusive for the container relocation operations in the yard, which are the configuration variables and the movement variables. The second group includes those exclusive for the stowage planning, and the third group includes the variables that integrate the two problems. All these variables are detailed next.

\begin{itemize}
\item  \textit{Yard configuration variables}: They indicate where the containers are located at any given time period $t$.
\end{itemize}

$b_{ijnt}= \left\{
\begin{array}{ll}
 1: & \text{if position $(i,j)$ is occupied by container $n$ in time period $t$,}  \\ 
 0: & \text{otherwise.}   \\
\end{array}
\right.$ 
\newline

$v_{nt}= \left\{
\begin{array}{ll}
 1: &\multirow{2}{*}{\parbox{9.0cm}{if container $n$ has been retrieved in time period $t'$, $t' = t+1$, with $t' \in \{1,...,t-1 \}$,}} \\
    & \\
 0: & \text{otherwise.}   \\
\end{array}
\right.$

\begin{itemize}
\item \textit{Movement variables}: Account for movements, either within the stack or from the stack to the outside.
\end{itemize}
  
$x_{ijklnt}= \left\{
\begin{array}{ll}
 1: & \multirow{2}{*}{\parbox{9.0cm}{if container $n$ is relocated from position $(i,j)$ to $(k,l)$ in time period $t$,}} \\
    &  \\
 0: & \text{otherwise.}   \\
\end{array}
\right.$ 
\newline

$y_{ijnt}= \left\{
\begin{array}{ll}
 1: & \text{if container $n$ is retrieved from position $(i,j)$ in time period $t$,} \\
 0: & \text{otherwise.}   \\
\end{array}
\right.$
\newline

\begin{itemize}
\item  \textit{Variables for the ship loading}: They integrate the movement between the yard and the ship.
\end{itemize}

$z_{ntrc}= \left\{
\begin{array}{ll}
 1: & \multirow{2}{*}{\parbox{9.0cm}{if container $n$ occupies the position $(r,c)$ in period $t$ into the ship,}}  \\
    &  \\
 0: & \text{otherwise.}  \\
\end{array}
\right.$ 

$q_{odrc}= \left\{
\begin{array}{ll}
 1: & \multirow{2}{*}{\parbox{9.0cm}{if a container with final destination port $d$ is relocated to position $(r,c)$ in port $o$,}} \\
    & \\
 0: & \text{otherwise.}   \\
\end{array}
\right.$

\begin{itemize}
\item  \textit{Ship variables}: They indicate the movement on the ship.
\end{itemize}

$w_{odarc}= \left\{
\begin{array}{ll}
 1: & \multirow{2}{*}{\parbox{9.0cm}{if there is a container in position $(r,c)$, loaded on board in port $o$, with final destination $d$, and unloaded in port $a$,}} \\
    &   \\
 0: & \text{otherwise.}   \\
\end{array}
\right.$ 
\newline

$u_{orc}= \left\{
\begin{array}{ll}
 1: &  \multirow{2}{*}{\parbox{9.0cm}{if, upon departing from port $o$, the position $(r, c)$ is occupied by a container,}} \\
    &  \\
 0: & \text{otherwise.}   \\
\end{array}
\right.$
\newline

As in \cite{Avriel1998}, the possible values for the index $a$ are  $\{o + 1, o + 2, ..., d - 1, d\}$. If $a = d$, the container is unloaded at its final destination, whereas if $a < d$, the container is relocated, that is, it is unloaded in port $a$ and reloaded again in the same port. It will then be denoted $w_{ada'r'c'}$, where $(r', c')$ may be a different position in the ship.

\subsection*{Objective function}

The model's objective function is to minimize the total number of relocation movements, which include the number of relocations in the yards (represented by the variable $x_{ijklnt}$), and the number of relocations at the ship (represented by the variable $w_{oadrc}$).

\begin{equation}
\begin{split}
\label{eq_obj}
Min \sum_{o=1}^{P-1}\left( \sum_{i=1}^{W_{o}}\sum_{j=1}^{H_{o}} \sum_{k=1}^{W_{o}} \sum_{l=1}^{H_{o}}\sum_{n = 1 | \atop n \in \varOmega_{o}}^{N_{o}} \sum_{t=1}^{N_{o}}x_{ijklnt} \ + \sum_{a=o+1}^{d-1}\sum_{d=o+1}^{P}\sum_{r=1}^{R}\sum_{c=1}^{C} w_{oadrc}\right)  
\end{split}
\end{equation}

\subsection*{Constraints}

\begin{itemize}
\item[I:] \textit{Yard constraints}
\end{itemize}

\begin{equation} \label{eq1}
\begin{aligned}
\sum_{i = 1}^{W_{o}} \sum_{j = 1}^{H_o} b_{ijnt} + v_{nt} =1 \\ \text{\footnotesize $n \in \varOmega_{o}, \ t\in T_{o}, \ o\in \{1,...,P-1\} $} 
\end{aligned}
\end{equation}

\begin{equation} \label{eq2}
\begin{aligned}
\sum_{n \in \varOmega_{o}}^{N_{o}} b_{ijnt} \leq 1 \\ \text{\footnotesize $i \in W_{o}, \ j \in H_{o}, \ t\in T_{o}, \ o\in \{1,...,P-1\}$}
\end{aligned}
\end{equation}

\begin{equation} \label{eq3}
\begin{split}
\sum_{n = 1 | \atop n \in \varOmega_{o}}^{N_{o}} b_{ijnt} \geq  \sum_{n = 1 | \atop n \in \varOmega_{o}}^{N_{o}} b_{i(j+1)nt} \\ \text{\footnotesize $i \in W_{o}, \ t\in T_{o}, \ j\in \{1,..., H_{o}-1\}, \ o\in \{1,...,P-1\}$}
\end{split}
\end{equation}

\begin{equation} \label{eq4}
\begin{split}
b_{ijnt} =  \sum_{k = 1}^{W_{o}} \sum_{l=1}^{H_{o}} x_{klijn(t-1)} - \sum_{k=1}^{W_{o}} \sum_{l=1}^{H_{o}} x_{ijkln(t-1)} + ~b_{ijn(t-1)} - y_{ijn(t-1)} \\ \text{\footnotesize $n \in  \varOmega_{o}, \ i \in  W_{o}, \ j \in H_{o}, \ t \in \{2,...,T_{o}\}, \ o \in \{1,...,P-1\}$}
\end{split}
\end{equation}

\begin{equation} \label{eq5}
\begin{split}
v_{nt} = \sum_{i=1}^{W_{o}}\sum_{j=1}^{H_{o}}\sum_{t'=1}^{t-1} y_{ijnt'} \\ \text{\footnotesize $n \in \varOmega_{o}, \ t \in \{2,...,T_{o}\}, \ o \in \{1,...,P-1\}$}
\end{split}
\end{equation}

\begin{equation} \label{eq6}
\begin{split}
M \left( 1 - \sum_{n = 1 | \atop n \in \varOmega_{o}}^{N_{o}} x_{ijklnt} \right)  \geq \sum_{n = 1 | \atop n \in \varOmega_{o}}^{N_{o}} \sum_{j'=j+1}^{H_{o}}\sum_{l'=l+1}^{H_{o}} x_{ij'kl'nt} \\ \text{\footnotesize $i, k \in W_{o}, \ j, l \in H_{o}, \ t \in T_{o}, \ o \in \{1,...,P-1\}$}
\end{split}
\end{equation}

\begin{equation} \label{eq7}
\begin{split}
M \left( 1 - \sum_{j=1}^{H_{o}}b_{ijnt} \right)  \geq \sum_{i' = 1 | \atop i' \neq i}^{W_{o}} \sum_{j=1}^{H_{o}} \sum_{k=1}^{W_{o}} \sum_{l=1}^{H_{o}} \sum_{n = 1 | \atop n \in \varOmega_{o}}^{N_{o}} x_{i'jklnt} \\  \text{\footnotesize $n \in \varOmega_{o}, \ i \in W_{o}, \ t \in T_{o}, \ o \in \{1,...,P-1\}$}
\end{split}
\end{equation}

\begin{equation} \label{eq8}
\begin{split} 
x_{ijilnt} = 0 \\ \text{\footnotesize $n \in  \varOmega_{o}, \ i \in W_{o}, \ j,l \in  H_{o}, \ t \in T_{o}, \ o \in \{1,...,P-1\}$}
\end{split}
\end{equation}

\begin{equation} \label{eq9}
\begin{split}
\sum_{k = 1}^{W_{o}} \sum_{l = 1}^{H_o}\sum_{n = 1 | \atop n \in \varOmega_{o}}^{N_{o}} x_{i(j+1)klnt} \ -   \sum_{n =1 | \atop n \in \varOmega_{o}}^{N_{o}} b_{i(j+1)nt} \ + 1 \geq  \sum_{k = 1}^{W_{o}} \sum_{l = 1}^{H_o}\sum_{n = 1 | \atop n \in \varOmega_{o}}^{N_{o}} x_{ijklnt} \ + \sum_{n =1 | \atop n \in \varOmega_{o}}^{N_{o}} y_{ijnt} \\ \text{\footnotesize  $i \in W_{o}, \ j \in \{1,...,H_{o}-1\}, \ t \in T_{o}, \ o \in \{1,...,P-1\}$}
\end{split}
\end{equation}

\begin{itemize}
\item[II:] \textit{Integration constraints}
\end{itemize}

\begin{equation} \label{eq10}
\begin{split}
 \sum_{n = 1 | \atop n \in \varOmega_{o}}^{N_{o}} v_{nt} = t \\ \text{\footnotesize $t\in T_{o}, \ o \in \{1,...,P-1\}$}
 \end{split}
\end{equation}

 \begin{equation} \label{eq11}
\begin{split}
\ v_{n(t+1)} \geq v_{nt} \\ \text{\footnotesize $n \in \varOmega_{o}, \ t \in \{1,...,T_{o}-1\}, \ o \in \{1,...,P-1\}$}
\end{split}
\end{equation}

\begin{equation} \label{eq12}
\begin{split}
\sum_{r = 1}^{R}\sum_{c = 1}^{C} z_{ntrc} = v_{n(t+1)} \\ \text{\footnotesize $n \in \varOmega_{o}, \ t\in \{1,...,T_{o}-1\}, o\in \{1,...,P-1\}$}
\end{split}
\end{equation}

\begin{equation} \label{eq13}
\begin{split}
\sum_{o'=1}^{o-1} \sum_{d=o+1}^{P} \sum_{a=o+1}^{d} w_{o'darc} \ + \sum_{d=o+1}^{P} q_{odrc} \ +  \sum_{n =1 | \atop n \in \varOmega_{o}}^{N_{o}} z_{ntrc}  \leq  1 \\  \text{\footnotesize $r\in R, \ c \in C\ \ t \in T_{o}, \ o \in \{1,...,P-1\}$}
\end{split}
\end{equation}

\begin{equation} \label{eq14} 
\begin{split}
z_{n(t+1)rc} \geq  z_{ntrc} \\ \text{\footnotesize $n \in \varOmega_{o}, \ r \in R, \ c  \in C, t \in \{1,...,T_{o}-1\}, o\in \{1,...,P-1\}$}
\end{split}
\end{equation}

\begin{equation} \label{eq15}
\begin{split}
\sum_{a=o+1}^{d} w_{odarc} =  \sum_{ n = 1 | \atop n \in \phi_{od}^{n}}^{N_{o}} z_{nT_{o}rc} + q_{odrc} \\ \text{\footnotesize $r\in R, \ c \in C, \ d \in \{o+1,...,P\}, o \in \{1,...,P-1\}$}
\end{split}
\end{equation}

\begin{equation} \label{eq16}
\begin{split}
\sum_{r=1}^{R}\sum_{c=1}^{C} z_{nT_{o}rc} = 1 \\ \text{\footnotesize $n \in \varOmega_{o}, \ o \in \{1,...,P-1\}$}
\end{split}
\end{equation}

\begin{equation} \label{eq17}
\begin{split}
\sum_{n = 1 | \atop n \in \varOmega_{o}}^{N_{o}} z_{ntrc} + \sum_{o'=1}^{o-1}\sum_{d=o+1}^{P}\sum_{a=o+1}^{d} w_{o'darc} ~+ \sum_{d=o+1}^{P}q_{odrc} \geq \sum_{n =1 | \atop n \in \varOmega_{o}}^{N_{o}} z_{nt(r+1)c} \\ \text{\footnotesize $ t \in \{1,..,T_{o}\}, \ r \in \{1,...,R-1\}, \ c \in C, \ o \in \{1,...,P-1\}$}
\end{split}
\end{equation}

\begin{equation} \label{eq18}
\begin{split}
\sum_{o'=1}^{o-1}\sum_{c=1}^{C}\sum_{r=1}^{R} w_{o'dorc} = \sum_{r=1}^{R}\sum_{c=1}^{C}q_{odrc} \\  \text{\footnotesize $d \in \{o+1,...,P\}, \ o \in \{1,...,P-1\}$}
\end{split}
\end{equation}

\begin{itemize}
\item[III:] \textit{Ship constraints}
\end{itemize}

\begin{equation} \label{eq19}
\begin{split}
\sum_{a=o+1}^{d}\sum_{r=1}^{R}\sum_{c=1}^{C} w_{odarc} \ - \sum_{m=1}^{o-1}\sum_{r=1}^{R}\sum_{c=1}^{C} w_{mdorc} = F_{od}  \\  \text{\footnotesize $o \in \{1,...,P-1\}, \ d \in \{o+1,...,P\}$}
\end{split}
\end{equation}

\begin{equation} \label{eq20}
\begin{split}
\sum_{m=1}^{o}\sum_{d=o+1}^{P}\sum_{a=o+1}^{d} w_{mdarc} = u_{orc} \\ \text{\footnotesize $o \in \{1,...,P-1\}, r \in R, \ c \in C$}
\end{split}
\end{equation}

\begin{equation} \label{eq21}
\begin{split}
u_{orc} - u_{o(r+1)c} \geq 0  \\ 
\text{\footnotesize $o \in \{1,...,P-1\}, \ r \in \{1,..., R-1\}, \ c \in C$}
\end{split}
\end{equation}

\begin{equation} \label{eq22}
\begin{split}
\sum_{o = 1}^{d-1} \sum_{e = d}^{P} w_{oedrc} + \sum_{o = 1}^{d-1} \sum_{e = d+1}^{P}\sum_{a = d+1}^{e} w_{oea(r+1)c}   \leq 1 \\  \text{\footnotesize $d \in \{2,...,P\}, \ r \in \{1,...,R-1\}, c \in C$}
\end{split}
\end{equation}

\begin{equation} \label{eq:stability}
\begin{split}
\sum_{r= \left\lceil \theta_{o} /\ C\right\rceil + 1}^{R} \sum_{c=1}^{C}  u_{orc} = 0  \\ \text{\footnotesize $o \in \{1,...,P-1\}$}
\end{split}
\end{equation}

The first group of constraints, (\ref{eq1}) to (\ref{eq9}), are those for the container removal process from the yard in the ports where the ship is moored.

Constraints (\ref{eq1}) ensure that in each time period, each container must either be within the bay or in the outside region. Constraints (\ref{eq2}) ensure that in each time period, each position $(i,j)$ must be occupied by at most one container. Constraints (\ref{eq3}) ensure that there are no `holes' in the stacking area by restricting that if there is a container in the position $(i,j+1)$, then the lower position $(i,j)$ must also be occupied. Constraints (\ref{eq4}) are the flow balancing constraint between the configuration and movement variables. They link the yard's layout in time period $t$ with its layout in period $t+1$ through the shifts and removals executed. Constraints (\ref{eq5}) establish that all containers shall be removed from the yard. Constraints (\ref{eq6}) ensure the LIFO policy. This constraint is necessary because more than one relocation is performed per time period. This means that if in time period $t$, container $n_{1}$ is under container $n_{2}$ and container $n_{1}$ is moved, then in period $t+1$ container $n_{1}$ cannot be below container $n_{2}$. Constraints (\ref{eq7}) ensure that only containers found above the target container are allowed to be relocated.  Constraints (\ref{eq8}) ensure that no container can be relocated to a position in the same stack where it currently is. Constraints (\ref{eq9}) ensure that a container in position $(i,j)$ can only be moved if the container at the position $(i,j+1)$ is relocated. If the container located at the position $(i,j+1)$ is not relocated, then we have $b_{i(j+1)nt} = 1$ and $x_{i(j+1)klnt} = 0$. This forces the left hand side of the equation to become equal to zero. Consequently, the right hand side of the equation must also be zero. Then, no relocation or retrieval is allowed for the container at position $(i,j)$.

The second group of constraints, (\ref{eq10}) to (\ref{eq18}), are those responsible for the integration between the ports' yards and the ship. These constraints ensure the correct ship loading. Constraints (\ref{eq10}) act as a counter for the remaining containers for a period of time. That is, in time period $t=1$, one container must have left, in time period $t=2$, two containers must have left, and so on until all containers $N_{o}$ have left in time period $T_{o}$. 
Constraints (\ref{eq11}) ensure that one container is retrieved from the yard per period of time. When container $n$ is retrieved from the yard, variable $v_{nt}$ becomes equal to 1 and continues being 1 over the following time periods. Constraints (\ref{eq12}) ensure that container $n$ is loaded on the ship in time period $t$. Constraints (\ref{eq13}) ensure that a position $(r,c)$ on the ship can only be occupied by one container, whether it is a container that was loaded in the current port (port $o$), in a previous port (port $o-1$) or a container that was already on board the ship and is being shifted in port $o$. Constraints (\ref{eq14}) certify that container $n$, after being loaded, does not change its position while the ship is still moored at the same port. Constraints (\ref{eq15}) ensure that if there is a container in the ship's position $(r,c)$, it must be a container that has just been loaded, or a shifted container. Constraints (\ref{eq16}) ensure that all the $N_{o}$ containers from yard $o$ have already been loaded onto the ship. Constraints (\ref{eq17}) ensure that during the ship's loading process, no container is allocated in a floating position or occupies the position of a container that was already in the ship or was shifted. Constraints (\ref{eq18}) count the total number of containers that have been shifted in port $o$.

The third set of constraints, (\ref{eq19}) to (\ref{eq:stability}), are exclusive for the ship. Constraints (\ref{eq19}) are related to the container's conservation flow. It indicates that the total number of containers to be shipped in port $o$ to destination port $d$ should be equal to the number of containers that were loaded in ports $p \in \{1,...,P-1\}$ subtracting the containers that were unloaded in ports $p \in \{o+1,...,P\}$. Constraint \ref{eq20} ensures that each slot $(r,c)$ is always occupied by at most one container. Constraints (\ref{eq21}) are required to ensure that there are containers underneath the container that occupies the position $(r,c)$. Constraints (\ref{eq22}) are responsible for defining how a container can be unloaded from the ship in port $d$ by imposing that if a container occupies the position $(r,c)$, then it will be unloaded in port $d$, if there is no container in the position $(r+1,c)$ above it. Constraints (\ref{eq:stability}) maintain the ship's stability by imposing a cargo height limit. This height limit is the ceiling for dividing the number of containers on board the ship by the number of stacks. Thus, the containers are equally distributed in the ship.

\section{Solution Approach} \label{sec:solution-approach}


The solution approach proposed here consists of generating feasible initial solutions for MPSP-CRP and inserts them into the optimization solver along with the mathematical model. This procedure is known as warm-starts. The advantages of using warm-starts is that it allows the solver to eliminate portions of the search space, and thus may result in smaller search trees, possibly leading to better solutions more quickly compared to a search initiated by looking for feasible solutions.

To define a feasible initial solution quickly, two constructive heuristic rules have been developed based on the problem structure. The first heuristic called HR1 is inspired by the heuristic proposed by \cite{Caserta2012} and it was adapted to solve the integrated problem. The idea is to determine where a relocating container should be placed based on a computation of a stack score. The second heuristic (HR2) randomly chooses where to relocate a container. Both heuristics are explained in detail hereafter.

\subsection{Heuristic HR1}

Given an instance $I$, let $L$ be the set of the containers currently located on the top of each stack in the yard belonging to port $o$, where $o \in \{1,...,P-1\}$; i.e. $L$ is the set of containers that could be removed from the yard in port $o$ without relocating other containers. For example, in Figure \ref{Descricao_Problema_Rota}, $L=\{1, 14, 15, 16\}$.
To determine which container in $L$ should be retrieved first, it is verified which one had the most distant port of destination $d$. This information is available in tuple $\phi_{od}^{n}$. Supposing that the respective destination port of containers in set $L$ is $P = \{4,4,5,4\}$, respectively. The heuristic will choose container $\hat{n}= 15$ to be retrieved first and load it onto the ship, since its destination is the farthest one. In case of a tie in this step, the heuristic randomly chooses the first container on the list that has the most distant destination.

Notice that since we are only removing the containers that are on the top of each stack, no relocations are performed at the yard.
To decide in which position in the ship this container ($\hat{n}=15$) should be placed, a stack store is calculated based on the nearest destination port of the containers that are already on board in the ship. To do this, let us define $min(c)$ as the nearest destination port of the containers in stack $c$ in the ship, with $c \in \{C\}$. For example, in Figure \ref{Descricao_Problema_Descarregamento}, this score would be $min(1)=4$, $min(2)=4$, $min(3)=5$ and $min(4)=3$, so, $min(c)=\{4,4,5,3\}$, since container number 2 in stack $c=3$ was unloaded in Port 2 releasing container 5 on the top. In case of an empty stack, the calculation set $min(c) = P+1$, and in case of full a stack, the calculation is set at $min(c) = 0$.

If there is a stack where $min(c)$ is still greater or equal than the port of destination of container $\hat{n}$ ($min(c) \geq \hat{n}$), then check if this position respects the maximum height allowed for this port (constraints (\ref{eq:stability})), if it does, choose a stack where $min(c)$ is minimum. The aim in this step is not to place the incoming container above another one that is going to be unloaded from the ship in a previous port, avoiding additional relocation at the next ports. If there is no stack satisfying the first condition, then choose the stack where $min(c)$ is maximum as container $\hat{n}$ will cause a new forced relocation. The reason for this is to postpone the relocation to the maximum. Briefly, the choice of stack $c^{*}$ where container $\hat{n}$ should be allocated in the ship is made as follows:

\begin{equation} \label{eq:find-c-HR1}
    \begin{split}
    \text{\footnotesize $c^{*}=$} \left\{
    \begin{array}{ll}
     \text{\footnotesize $\argmin_{ c \in \{C\}} \ \{min(c) : min(c) \geq \hat{n} \}, \ \text{if}$} \ \text{\footnotesize $\exists \ c: min(c) \geq \hat{n};$} \\
     \text{\footnotesize $\argmax_{c \in \{C\}} \{min(c) \}, ~ \text{otherwise.}$}   \\
    \end{array}
    \right. 
    \end{split}
\end{equation}








After loading all containers from port $p$, if there is a container above the maximum height allowed, then a new position is chosen for it following the steps described previously. When loading is finished in port $p$, the ship goes to the next port $p+1$. It is worth remembering that from the second port to port $P-1$, the unloading of all containers destined to port $p$ is performed first, and then the loading of the yard containers is started

Therefore to unload, the heuristic rule verifies if the container with a destination to port $p$ is at the top of its stack. If it is, it is unloaded. If it is not, define the relocation set $Rel$ with the blocking containers above the target container according to the LIFO policy. Since there is no relocation allowed inside the ship, these containers are unloaded to a reserved area in the yard and then loaded again after all the unloading is performed. Algorithm \ref{algorithm:HR1} presents the heuristic steps as a pseudocode.

\begin{algorithm}[H] \label{algorithm:HR1}
\small
\begin{algorithmic}[1]
\linespread{0.90}\selectfont
\State $initialize(\text{problem instance})$ 
\State $p \gets 1$, $move\_counter \gets 0$
\For{$p = \{1,...,P\}$}
\If{$p = 1$}
\While{\text{Yard($p$)} $\neq \emptyset$}
\State define set $L$ and set $D$ 
\State $\hat{n} = max(D)$ 
\For{$c = \{1,...,C\}$}
\If{stack $c$ is empty}
\State $min(c) = P+1$
\EndIf
\If{stack $c$ is full}
\State $min(c) = 0$
\EndIf
\If{stack $c$ is not empty or full}
\State $min(c)=min\{k:$ port of destination $k$ in stack $c$\}
\EndIf
\EndFor
\State find $c^{*}$ using Eq. (\ref{eq:find-c-HR1})
\State \text{Yard($p$)} $\gets$ \text{Yard($p$)}$\backslash \{ \hat{n} \}$ 
\State put $\hat{n}$ in stack $c^{*}$ 
\EndWhile
\EndIf
\algstore{myalg}
\caption{Heuristic HR1}
\end{algorithmic}
\label{Simulation}
\end{algorithm}

\begin{algorithm}[t]                     
\begin{algorithmic}[1]                   
\linespread{0.90}\selectfont
\algrestore{myalg}
\State $p\gets p+1$
\If{$1 < p < P$}
\State define set $Unload$ \{set of containers to be unloaded in $p$\} and set $Temp\_Yard$ \{Yard to allocate relocated containers\}
\While{$Unload$ $\neq \emptyset$} 
\State set $\hat{n}$ as the first element in $Unload$ 
\If{$\hat{n}$ is at the top of its stack}
\State $Unload$ $\gets$ $Unload\backslash \{ \hat{n} \}$ 
\Else
\State define set $Rel$ 
\State set $\hat{r}$ as the first element in $Rel$ 
\State move $\hat{r}$ to $Temp\_Yard$ using the same procedure described when $p=1$
\State $Unload$ $\gets$ $Unload\backslash \{ \hat{r} \}$ 
\State $move\_counter\gets move\_counter + 1$
\EndIf
\EndWhile
\State move containers from $Temp\_Yard$ back to the ship using the same procedure described when $p=1$
\State repeat loading step as described when $p=1$
\If{Eq. (\ref{eq:stability}) is violated}
\State find new $c^{*}$ using Eq. (\ref{eq:find-c-HR1})
\State $move\_counter\gets move\_counter + 1$
\EndIf
\State $p\gets p+1$
\EndIf
\If{$p = P$}
\State perform only unloading procedure
\EndIf
\EndFor
\end{algorithmic}
\label{Simulation}
\end{algorithm}

\subsection{Heuristic HR2}

What differentiates the heuristic rule HR2 from the HR1 is the way of choosing a position for the container that is being loaded into the ship when there are no positions that will not cause future relocations. While HR1 uses Equation \ref{eq:find-c-HR1}, which aims to delay future relocations to the maximum, HR2 randomly chooses an available position. Thus, the choice of the $c^{*}$ column in which the $\hat{n}$ container is going to be allocated is made as follows by HR2:

\begin{equation} \label{eq:find-c-HR2}
    \begin{split}
    \text{\footnotesize $c^{*}=$} \left\{
    \begin{array}{ll}
     \text{\footnotesize $\argmin_{ c \in \{C\}} \ \{min(c) : min(c) \geq \hat{n} \}, \ \text{if}$} \ \text{\footnotesize $\exists \ c: min(c) \geq \hat{n};$} \\
     \text{\footnotesize $\datasample_{c \in \{C\}} \{min(c) \}, ~ \text{otherwise.}$}   \\
    \end{array}
    \right. 
    \end{split}
\end{equation}

The $\datasample(min(c),1)$ function returns a observation sampled uniformly randomly in $min(c)$. In Algorithm \ref{algorithm:HR1}, the modification would be in lines 19 and 42, which would be replaced by: ``Find $c^{*}$'' using  Eq. (\ref{eq:find-c-HR2}).

The result of the HR1 is set as the upper bound for HR2. To ensure that at any time HR2 finds a better solution than HR1, the HR2 runs 10,000 times. If a better solution shows up, then this solution is kept until a better one is obtained. The motivation behind the modification proposed by the HR2 is based on the hypothesis that the smaller the gap between the initial solution generated by the heuristic and the value of the optimal solution, the greater the impact on the warm-start initialization, and on the computational time required to find the optimal solution.

\subsection{Warm-start Initialization}

To use the result of the heuristic rules for warm-starting the exact method represented by Cplex\texttrademark{} and Gurobi\texttrademark{}, while the heuristic is running, the information related to the movement of the containers in each step is saved in the matrix form. This information is then translated into the binary variables of the MPSP-CRP to provide a warm-start as a parameter for the model. It is worth remembering that the proposed heuristic always generate feasible solutions for the model. Next, the computational results are presented.

\section{Computational Results} \label{sec:computational-results}

In this section, the computational results of the MPSP-CRP without a warm-start initialization are presented first, followed by those from the the hierarchical approach, and the HR1 and HR2 heuristics. Finally, the results from the MPSP-CRP with a warm-start initialization are presented. 

\subsection{MPSP-CRP}

The MPSP-CRP model was implemented in Python and solved in a computer with an AMD Ryzen 7 1700 processor, 3 GHz, 64 GB RAM, using 16 threads.
 
Tables \ref{tabInstancias1} and \ref{tabInstancias2} describe the randomly generated test instances. For each instance $I$ the following is described: number of yard tiers ($H$), number of yard stacks ($W$), number of ship tiers ($R$), number of ship stacks ($C$), number of ports traveled by the ship ($P$), total number of containers ($N$) with which the ship handles along its entire route, and finally, the total number of variables and constraints that this instance produces in the exact model. The yards of all instances have an occupancy rate randomly generated between 65\% to 85\%. For each instance of Table \ref{tabInstancias1}, the ship created has enough space to load all containers during its route.

Note that the instances in Table \ref{tabInstancias1} are the same in Table \ref{tabInstancias2}, but with a bigger number of stacks on the ship. The impact of having more space at the ship to handle the containers can be observed by doing this. Note also that increasing the ship's number of stacks significantly increases the size of the problem in terms of the number of variables and constraints.

Instances $12A$ and $12B$ from Tables \ref{tabInstancias1} and \ref{tabInstancias2} appear with a "$\ast$" in the number of variable and constraint columns because the available memory was not sufficient to finish writing the models, besides using procedures to save memory such as those recommended by the solvers. All instances are available in the repository https://github.com/JunqueiraCatarina/Instances.

All instances in Tables \ref{tabInstancias1} and \ref{tabInstancias2} were solved in two different solvers that are the best academically and commercially known to deal with linear and integer programming models: Cplex\texttrademark{} 12.8 and Gurobi\texttrademark{} 8.0. The motivation for using these two solvers is to compare their results and check if one of them could deal more efficiently with the proposed model as it was clearly a difficult problem to solve.

\begin{table}[h!]
    \caption{Test instances. \label{tabInstancias1}}
    \begin{center}
    \scalebox{0.80}{
    \begin{tabular}{|c|cccc|cccc|c|c|c|c|}
      \hline
      \multirow{2}{*}{$I$} & & \multicolumn{2}{c}{Yard} & & & \multicolumn{2}{c}{Ship} & & \multirow{2}{*}{$P$} &  \multirow{2}{*}{\makecell{$N$}} & \multirow{2}{*}{\makecell{Number of\\ Variables}} &  \multirow{2}{*}{\makecell{Number of\\ Restrictions}} \\
     \cline{3-4} \cline{7-8}
        & & $H$ & $W$ & & & $R$ & $C$ & & & & & \\
      \hline
      $1A$ & & 4 & 3 & & & 4 & 5  & & 4   &  30  & 57.080   & 30.373 \\
      $2A$ & & 4 & 3 & & & 4 & 6  & & 6   &  46  & 83.924   & 45.864 \\  
      $3A$ & & 4 & 3 & & & 5 & 9  & & 10  &  92  & 213.155  & 123.772 \\
      $4A$ & & 4 & 3 & & & 6 & 10 & & 15  &  137 & 351.493  & 204.102 \\
      \hline
      $5A$ & & 4 & 6 & & & 5 & 5  & & 4   &  58 & 757.955    & 223.701 \\
      $6A$ & & 4 & 6 & & & 6 & 6  & & 6   &  91 & 1.125.877  & 338.455 \\  
      $7A$ & & 4 & 6 & & & 7 & 9  & & 10  & 162 & 2.053.070  & 665.730 \\
      $8A$ & & 4 & 6 & & & 7 & 10 & & 15  & 261 & 3.773.715  & 1.407.996 \\
      \hline
      $9A$  & & 6 & 10 & & & 6 & 14 & & 4   &  131 & 21.861.321  & 3.340.081  \\
      $10A$ & & 6 & 10 & & & 6 & 19 & & 6   &  212 & 34,582,630  & 5,531,912  \\  
      $11A$ & & 6 & 10 & & & 6 & 30 & & 10  &  409 & 73,054,437  & 12,664,333 \\
      $12A$ & & 6 & 10 & & & 6 & 44 & & 15  &  617 & $\ast$      & $\ast$     \\
      \hline
    \end{tabular}}
    \end{center}
\end{table}

\begin{table}[H]
    \caption{Test instances. \label{tabInstancias2}}
    \begin{center}
    \scalebox{0.80}{
    \begin{tabular}{|c|cccc|cccc|c|c|c|c|}
  \hline
  \multirow{2}{*}{$I$} & & \multicolumn{2}{c}{Yard} & & & \multicolumn{2}{c}{Ship} & & \multirow{2}{*}{$P$} &  \multirow{2}{*}{\makecell{$N$}} & \multirow{2}{*}{\makecell{Number of\\ Variables}} &  \multirow{2}{*}{\makecell{Number of\\ Restrictions}} \\
    \cline{3-4} \cline{7-8}
   & & $H$ & $W$ & & & $R$ & $C$ & & & & & \\
  \hline
  $1B$ & & 4 & 3 & & & 4 & 10 & &  4   &  30 & 63,460  & 37,095 \\
  $2B$ & & 4 & 3 & & & 4 & 12 & &  6   &  46 & 95,516  & 57,627 \\  
  $3B$ & & 4 & 3 & & & 4 & 18 & &  10  &  92 & 265,760 & 172,915 \\
  $4B$ & & 4 & 3 & & & 4 & 20 & &  15  &  137& 473,653 & 300,916 \\
  \hline
  $5B$ & & 4 & 6 & & & 5 & 18 & &  4   & 58  & 809,660 & 277,261 \\
  $6B$ & & 4 & 6 & & & 6 & 16 & &  6   & 91  & 1,208,629 & 423,569 \\  
  $7B$ & & 4 & 6 & & & 7 & 20 & &  10  & 162 & 2,273,640 & 885,553 \\
  $8B$ & & 4 & 6 & & & 7 & 36 & &  15  & 261 & 4,478,055 & 2,072,990 \\
  \hline
  $9B$  & & 6 & 10 & & & 6 & 28 & &  4   & 131 & 22,345,497 & 3,833,009 \\
  $10B$ & & 6 & 10 & & & 6 & 38 & &  6   & 212 & 35,616,724 & 6,578,547 \\  
  $11B$ & & 6 & 10 & & & 6 & 60 & &  10  & 409 & 76,462,917 & 16,099,967 \\
  $12B$ & & 6 & 10 & & & 6 & 88 & &  15  & 617 & $\ast$     & $\ast$  \\
  \hline
  \end{tabular}}
  \end{center}
\end{table}

Given this, Tables \ref{tabResultadosSolverInstancias1} and \ref{tabResultadosSolverInstancias2} report the computational results from the solvers throughout the experiments under a time limit of 432,000 seconds (5 days). After that, the solver is stopped and the result reached until then is reported. The first column, $I$, shows the number of the instance, the $O.F.$ illustrates the value of the objective function, the $Time$ column shows the time in seconds that the respective solver took to solve the instance, and the $Gap$ column depicts the gap between the solution found and the best bound. Note that in the instances with a "$-$" symbol in $O.F.$ and $Gap$ columns are the ones where the solver reached the maximum computational time of 432,000 seconds without finding a feasible solution.

The last four instances of both Tables \ref{tabResultadosSolverInstancias1} and \ref{tabResultadosSolverInstancias2} are filled with a "$\ast$" symbol. This means that the 64GB of available RAM memory was not enough to solve these instances in the Branch and Bound procedure. This memory problem happens not because of the size of the problem itself, as shown in Tables \ref{tabInstancias1} and \ref{tabInstancias2}, but because of its combinatorial characteristic, whose list of active nodes in the Branch and Bound tree can consume large amounts of memory. This fact leads the solvers to terminate the optimization process with an error message. Moreover, the adjustment of memory emphasis parameters was not successful in these cases.

Based on the results analysis, it can be observed that Cplex\texttrademark{} and  Gurobi\texttrademark{} were able to find a solution for the same set of instances; however, Gurobi\texttrademark{} performed better in all instances where a solution was found. In all instances where the optimal solution was proved (gap of 0.00\%), Gurobi\texttrademark{} consumed less time than the Cplex\texttrademark{} to finish the optimization. Moreover, Gurobi\texttrademark{} found a solution with a lower gap in instance $4A$. In instance $7B$, Gurobi\texttrademark{} found the optimal solution in 16.98 hours; meanwhile, Cplex\texttrademark{} terminated the optimization because it reached the time limit of 5 days, and with a gap of 99.31\%.

When comparing the results in Table \ref{tabResultadosSolverInstancias1} with the those obtained in the instances where there is more handling space on the ship (Table \ref{tabResultadosSolverInstancias2}), can be observed that even the problem resulting in a larger number of variables and constraints, the solvers performed better: solutions were found in more instances and less time was used. This happens because it is easier to find solutions when there is more space on the ship, and less computational effort is required. In general, it is notable that the difficulty of the problem grows rapidly and that an instance with less than 200 containers, although quite distant from a real problem, already becomes prohibitive for the exact model until now.

 
 

 
 
 

\begin{table}[H]
\caption{Results: Cplex\texttrademark{} x Gurobi\texttrademark{}. \label{tabResultadosSolverInstancias1}}
\begin{center}
\scalebox{0.85}{
\begin{tabular}{|c|ccccc|ccccc|}
  \hline
  \multirow{2}{*}{$I$} & & \multicolumn{3}{c}{Cplex\texttrademark{}}  & & & \multicolumn{3}{c}{Gurobi\texttrademark{}}  &    \\
 \cline{3-5} \cline{8-10}
   & & $O.F.$ & $Time$ & $Gap$ & &  & $O.F.$ & $Time$ & $Gap$ &    \\
  \hline
  $1A$  & & 0 & 10.94      & 0.00\%  & & & 0 & 6.28       & 0.00\% &  \\
  $2A$  & & 1 & 36.42      & 0.00\%  & & & 1 & 22.95      & 0.00\% &  \\  
  $3A$  & & 2 & 13,669.55  & 0.00\%  & & & 2 & 3,260.94   & 0.00\% &  \\
  $4A$  & & 16 & 432,010.98 & 93.75\% & & & 3 & 432,000.62 & 66.67\%&  \\
  \hline
  $5A$ & & 0 & 6,079.77   & 0.00\% & & & 0 & 3,892.21 & 0.00\% &  \\
  $6A$ & & 0 & 21,351.66  & 0.00\% & & & 0 & 5,883.63 & 0.00\% &  \\  
  $7A$ & & - & 432,014.87 & -      & & & - & 432,001.47 & -    &  \\
  $8A$ & & - & 432,005.72 & -      & & & - & 432,000.85 & -    &  \\
  \hline
  $9A$  & & $\ast$ & $\ast$ & $\ast$ & & & $\ast$ & $\ast$ & $\ast$ &  \\
  $10A$ & & $\ast$ & $\ast$ & $\ast$ & & & $\ast$ & $\ast$ & $\ast$ &  \\  
  $11A$ & & $\ast$ & $\ast$ & $\ast$ & & & $\ast$ & $\ast$ & $\ast$ &  \\
  $12A$ & & $\ast$ & $\ast$ & $\ast$ & & & $\ast$ & $\ast$ & $\ast$ &  \\
  \hline
  \end{tabular}}
\end{center}
\end{table}

\begin{table}[H]
\caption{Results: Cplex\texttrademark{} x Gurobi\texttrademark{}. \label{tabResultadosSolverInstancias2}}
\begin{center}
\scalebox{0.85}{
\begin{tabular}{|c|ccccc|ccccc|}
  \hline
  \multirow{2}{*}{$I$} & & \multicolumn{3}{c}{Cplex\texttrademark{}}  & & & \multicolumn{3}{c}{Gurobi\texttrademark{}}  &    \\
 \cline{3-5} \cline{8-10}
   & & $O.F.$ & $Time$ & $Gap$ & &  & $O.F.$ & $Time$ & $Gap$ &    \\
  \hline
  $1B$ & & 0 & 6.11      & 0.00\%  & & & 0 & 3,30     & 0.00\% &  \\
  $2B$ & & 0 & 12.69     & 0.00\%  & & & 0 & 7,25     & 0.00\% &  \\  
  $3B$ & & 0 & 630.81    & 0.00\%  & & & 0 & 273,07   & 0.00\% &  \\
  $4B$ & & 0 & 43,951.89 & 0.00\%  & & & 0 & 3.872,97 & 0.00\% &  \\
  \hline
  $5B$ & & 0    & 262.83     & 0.00\%  & &  & 0 & 552.66     & 0.00\% &  \\
  $6B$ & & 0    & 2,401.53   & 0.00\%  & &  & 0 & 1,567.67   & 0.00\% &  \\  
  $7B$ & & 1,146 & 432,002.89 & 99.31\% & &  & 11& 61,126.12  & 0.00\% &  \\
  $8B$ & & -    & 432,002.89 & -       & &  & - & 432,000.88 & - &   \\
  \hline
  $9B$  & & $\ast$ & $\ast$ & $\ast$ & & & $\ast$ & $\ast$ & $\ast$ &  \\
  $10B$ & & $\ast$ & $\ast$ & $\ast$ & & & $\ast$ & $\ast$ & $\ast$ &  \\  
  $11B$ & & $\ast$ & $\ast$ & $\ast$ & & & $\ast$ & $\ast$ & $\ast$ &  \\
  $12B$ & & $\ast$ & $\ast$ & $\ast$ & & & $\ast$ & $\ast$ & $\ast$ &  \\
  \hline
  \end{tabular}}
\end{center}
\end{table}


 \subsection{Hierarchical MPSP-CRP approach} 

For the hierarchical approach, the MPSP is solved first, which consists of minimizing the variable $w_{odarc}$, subject to constraints (\ref{eq19}) to (\ref{eq:stability}). Then, the solution obtained for the variables $w_{odarc}$ and $u_{orc}$ is fixed and inserted in the rest of the problem, which consists of removing the containers from the yards and loading then on to their positions in the ship. This second part will be called CRP $+$ Loading, and is represented by the minimization of the variable $x_{ijklnt}$, subject to constraints (\ref{eq1}) to (\ref{eq18}).

Tables \ref{tabResultadosHierarquicosInstancias1} and \ref{tabResultadosHierarquicosInstancias2} report separately the results of the two steps of the hierarchical solution approach for the proposed instances solved by the Cplex\texttrademark{}. Since the objective of this experiment is to demonstrate the difference between the two approaches (hierarchical and integrated), only one of the solvers needs to be used. The \textit{Total O.F.} column shows the sum of the solutions obtained by each part, and the last column, \textit{Best MPSP-CRP O.F.} brings the best solution obtained for the MPSP-CRP, to facilitate the comparison between the two approaches.

As expected, even though this approach could solve more instances in less computational time, the hierarchical solution did not lead to as good solutions as the optimal solution of the integrated problem.

It can be observed that for instances $3A$, $4A$ and $4B$, the hierarchical optimal solution was worse than the integrated solution. This happened because a stowage planning made without taking the yard into consideration can lead to more relocations in the yard.

\begin{table}[H]
\caption{Results of the hierarchical approach using Cplex\texttrademark{}. \label{tabResultadosHierarquicosInstancias1}}
\begin{center}
\scalebox{0.80}{
\begin{tabular}{|c|ccccc|ccccc|c|c|}
  \hline
  \multirow{2}{*}{$I$} & & \multicolumn{3}{c}{MPSP}  & & & \multicolumn{3}{c}{CRP $+$ Loading} & &  \multirow{2}{*}{Total O.F.} & Best MPSP- \\
 \cline{3-5} \cline{8-10}
   & & $O.F.$ & $Time$ & $Gap$ & &  & $O.F.$ & $Time$ & $Gap$ & & & CRP $O.F.$  \\
  \hline
  $1A$  & & 0 & 0.14     & 0.00\%  & & & 0 & 6.05   & 0.00\% & & 0 & 0 \\
  $2A$  & & 1 & 0.24     & 0.00\%  & & & 0 & 6.47   & 0.00\% & & 1 & 1 \\  
  $3A$  & & 2 & 2.47     & 0.00\%  & & & 1 & 48.52  & 0.00\% & & 3 & 2 \\
  $4A$  & & 1 & 5,701.60 & 0.00\%  & & & 8 & 588.76 & 0.00\% & & 9 & 3 \\
  \hline
  $5A$ & & 0  & 0.04      & 0.00\%  & & & 0 & 249.2     & 0.00\% & & 0 & 0 \\
  $6A$ & & 0  & 0.16      & 0.00\%  & & & 0 & 423.74    & 0.00\% & & 0 & 0 \\  
  $7A$ & & 12 & 43,227.55 & 11.67\% & & & 2 & 10,927.56 & 0.00\% & & 14 & - \\
  $8A$ & & 1  & 5,602.46  &  0.00\% & & & 2 & 36,398.90 & 0.00\% & & 3 & - \\
  \hline
  $9A$  & & 0      & 0.08   & 0.00\% & & & $\ast$ & $\ast$ & $\ast$ & & $\ast$ & $\ast$ \\
  $10A$ & & $\ast$ & $\ast$ & $\ast$ & & & $\ast$ & $\ast$ & $\ast$ & & $\ast$ & $\ast$ \\  
  $11A$ & & $\ast$ & $\ast$ & $\ast$ & & & $\ast$ & $\ast$ & $\ast$ & & $\ast$ & $\ast$ \\
  $12A$ & & $\ast$ & $\ast$ & $\ast$ & & & $\ast$ & $\ast$ & $\ast$ & & $\ast$ & $\ast$ \\
  \hline
  \end{tabular}}
\end{center}
\end{table}

\begin{table}[H]
\caption{Results of the hierarchical approach using Cplex\texttrademark{}. \label{tabResultadosHierarquicosInstancias2}}
\begin{center}
\scalebox{0.80}{
\begin{tabular}{|c|ccccc|ccccc|c|c|}
  \hline
  \multirow{2}{*}{$I$} & & \multicolumn{3}{c}{MPSP}  & & & \multicolumn{3}{c}{CRP $+$ Loading} & &  \multirow{2}{*}{Total $O.F.$} & Best MPSP-  \\
 \cline{3-5} \cline{8-11}
   & & $O.F.$ & $Time$ & $Gap$ & &  & $O.F.$ & $Time$ & $Gap$ & & & CRP $O.F.$ \\
  \hline
  $1B$  & & 0 & 0.10 & 0.00\%  & & & 0 & 6.16   & 0.00\% & & 0 & 0 \\
  $2B$  & & 0 & 0.13 & 0.00\%  & & & 0 & 10.59  & 0.00\% & & 0 & 0 \\  
  $3B$  & & 0 & 0.73 & 0.00\%  & & & 0 & 39.71  & 0.00\% & & 0 & 0 \\
  $4B$  & & 0 & 9.52 & 0.00\%  & & & 2 & 205.44 & 0.00\% & & 2 & 0 \\
  \hline
  $5B$ & & 0  & 0.04  & 0.00\% & & & 0 & 162.6   & 0.00\% & & 0 & 0 \\
  $6B$ & & 0  & 0.16  & 0.00\% & & & 0 & 189.24  & 0.00\% & & 0 & 0 \\  
  $7B$ & & 11 & 7.79  & 0.00\% & & & 0 & 832.62  & 0.00\% & & 11 & 11 \\
  $8B$ & & 1  & 15.68 & 0.00\% & & & 0 & 2,109.6 & 0.00\% & & 1 & - \\
  \hline
  $9B$  & & 0      & 0.05   & 0.00\% & & & $\ast$ & $\ast$ & $\ast$ & & $\ast$ & $\ast$ \\
  $10B$ & & $\ast$ & $\ast$ & $\ast$ & & & $\ast$ & $\ast$ & $\ast$ & & $\ast$ & $\ast$ \\  
  $11B$ & & $\ast$ & $\ast$ & $\ast$ & & & $\ast$ & $\ast$ & $\ast$ & & $\ast$ & $\ast$ \\
  $12B$ & & $\ast$ & $\ast$ & $\ast$ & & & $\ast$ & $\ast$ & $\ast$ & & $\ast$ & $\ast$ \\
  \hline
  \end{tabular}}
\end{center}
\end{table}

\subsection{Heuristics based on the problem}

Table \ref{tabResultadosHeuristica} reports the results from the heuristic previously described. The results are shown in terms of the objective function ($O.F.$), which is the total number of relocations, and $Time$, which is the total computational time in seconds required to solve each instance ($I$).

By analyzing Table \ref{tabResultadosHeuristica}, one can observe that heuristic HR1 was able to solve most of the instances in less than one second, while HR2 was slightly more time-consuming. Nevertheless, the latter had better results in the objective function for some instances, showing that postponing the relocation may not always be the best choice.

\begin{table}[H]
\caption{Results for the HR1 and HR2 heuristics. \label{tabResultadosHeuristica}}
\begin{center}
\begin{tabular}{cc}
\scalebox{0.80}{
\begin{tabular}{|c|cc|cc|}
  \hline
  \multirow{2}{*}{$I$} &  \multicolumn{2}{c|}{HR1}  & \multicolumn{2}{c|}{HR2}    \\
  \cline{2-3} \cline{4-5}
   & $O.F.$ & $Time$ & $O.F.$ & $Time$  \\
  \hline 
  $1A$ &  1 & 0.13  & 1  & 53.72 \\
  $2A$ &  4 & 0.14  & 4  & 111.35  \\  
  $3A$ & 14 & 0.18  & 8  & 282.55 \\
  $4A$ & 20 & 0.20  & 16 & 804.34   \\
  \hline
  $5A$ & 0  & 0.11  & 0  &  38.02 \\
  $6A$ & 5  & 0.16  & 5  & 230.18 \\  
  $7A$ & 66 & 0.16  & 59 & 840.44 \\
  $8A$ & 63 & 0.27  & 48 & 1,339.20 \\
  \hline
  $9A$  & 0  & 0.12 & 0  & 108.46 \\
  $10A$ & 12 & 0.22 & 11 & 882.72  \\  
  $11A$ & 33 & 0.35 & 26 & 2,625.30  \\
  $12A$ & 82 & 0.68 & 66 & 5,708.80  \\
  \hline
\end{tabular}}
&
\scalebox{0.80}{
\begin{tabular}{|c|cc|cc|}
  \hline
  \multirow{2}{*}{$I$} &  \multicolumn{2}{c|}{HR1}  & \multicolumn{2}{c|}{HR2}    \\
  \cline{2-3} \cline{4-5}
   & $O.F.$ & $Time$ & $O.F.$ & $Time$  \\
  \hline 
  $1B$ &  0  & 0.11  &  0   & 21.88 \\
  $2B$ &  0  & 0.12  &  0   & 18.86  \\  
  $3B$ &  7  & 0.17  &  7   & 432.74 \\
  $4B$ &  8  & 0.16  &  8   & 946.82   \\
  \hline
  $5B$ & 0  & 0.12   &  0  & 51.96 \\
  $6B$ & 0  & 0.12   &  0  & 42.14 \\  
  $7B$ & 44 & 0.26   & 40  & 1,079.40 \\
  $8B$ & 32 & 0.35   & 32  & 2,257.70 \\
  \hline
  $9B$  & 0  & 0.16 &   0     & 195.03 \\
  $10B$ & 0  & 0.26 &   0     & 227.87 \\  
  $11B$ & 0  & 0.60 &   0     & 345.05  \\
  $12B$ & 34 & 1.48 &   34    & 1,2841.0  \\
  \hline
\end{tabular}}
\end{tabular}
\end{center}
\end{table}

To demonstrate the efficiency of the heuristic procedures, they were also tested in two sets of larger instances. The first set contains instances up to 10,000 containers and all yards have the same dimensions. These new instances and their results are presented in Table \ref{tabInstanciasGrandes1}. For this set of instances, a maximum computational time of two hours (7,200 seconds) or 10,000 rounds was set for heuristic HR2. When any of the conditions is reached, the run is stopped and the best solution obtained until then is reported.

The second set shows real-world scaled instances, where the yards have sizes equivalent to a bay of the ports of Aarhus in Denmark; Antwerp in Belgium; Busan in Korea; Hamburg in Germany; Jebel Ali in Dubai; Rotterdam in Netherlands; Santos in Brazil; Shanghai in China; and Singapore, and the ship can load up to 30,000 containers simultaneously. These new instances and their results are presented in Table \ref{tabInstanciasGrandes2}. For this set of instances, a maximum computational time of 12 hours (43,200 seconds) or 10,000 rounds was set for the HR2 heuristic.

The HR1 heuristic was able to solve all proposed instances in reasonable computational time, even the real-world scaled ones. Note that instances up to 10,000 containers were solved in less than two minutes. Meanwhile heuristic rule HR2 reached the maximum computational time allowed in all instances. Despite this, the rule was able to improve the solution of some of the instances of Table \ref{tabInstanciasGrandes1}. The time may seem long, but it is still much less than the exact approach took for much smaller instances. In short, the heuristic HR1 turned out to be a more efficient solving methodology, since HR2 was not able to improve the solution of the real-world scaled instances.

\begin{table}[b!]
\caption{Large instances and computational results from the heuristics. \label{tabInstanciasGrandes1}}
\begin{center}
\scalebox{0.79}{
\begin{tabular}{|c|cccc|cccc|c|c|ccccccc|}
  \hline
   \multirow{2}{*}{$I$} & & \multicolumn{2}{c}{Yard} & & & \multicolumn{2}{c}{Ship} & & \multirow{2}{*}{$P$} &  \multirow{2}{*}{\makecell{$N$}} & & \multicolumn{2}{c}{HR1} & & \multicolumn{2}{c}{HR2} &  \\
 \cline{3-4} \cline{7-8} \cline{13-14} \cline{16-17}
   & & $H$ & $W$ & & & $R$ & $C$ & & & & & $O.F.$ &  $Time$ & & $O.F.$ & $Time$ & \\
  \hline
  $13A$ & & 10 & 100 & & & 28 & 73  & &  4   &  2,352  & & 164 & 5.56     & & 164  & 7,200 & \\
  $14A$ & & 10 & 100 & & & 33 & 91  & &  6   &  3,801  & & 872 & 16.65    & & 815  & 7,200 & \\  
  $15A$ & & 10 & 100 & & & 39 & 108 & &  10  &  6,517  & & 2,110 & 42.86  & & 2033 & 7,200 & \\
  $16A$ & & 10 & 100 & & & 49 & 133 & &  15  &  10,372 & & 4,512 & 99.25  & & 4512 & 7,200 & \\
  \hline
  $13B$ & & 10 & 100 & & & 28 & 146  & &  4   & 2,352 & & 45   & 11.02   & & 45    & 7,200 & \\
  $14B$ & & 10 & 100 & & & 33 & 182  & &  6   & 3,801 & & 234  & 21.31   & & 221   & 7,200 & \\  
  $15B$ & & 10 & 100 & & & 39 & 216  & &  10  & 6,517 & & 749  & 52.75   & & 749   & 7,200 & \\
  $16B$ & & 10 & 100 & & & 49 & 266  & &  15  & 10,372& & 1,525 & 107.18 & & 1,499 & 7,200 & \\
  \hline
\end{tabular}}
\end{center}
\end{table}

\begin{table}[b!]
\caption{Set of real-scale instances and results from the heuristics. \label{tabInstanciasGrandes2}}
\begin{center}
\scalebox{0.79}{
\begin{tabular}{|c|cccc|c|c|cccc|cccc|}
  \hline
   \multirow{2}{*}{$I$} & & \multicolumn{2}{c}{Ship} & & \multirow{2}{*}{$P$} &  \multirow{2}{*}{\makecell{$N$}} & & \multicolumn{2}{c}{HR1} & & & \multicolumn{2}{c}{HR2} &  \\
 \cline{3-4} \cline{9-10} \cline{13-14} 
        &  & $R$ & $C$ & &    &        & & $O.F.$ &  $Time$  & &  & $O.F.$ & $Time$ & \\
  \hline
  $17A$ &  & 45 & 551 & &  6  & 30,306 & & 8,384  & 3,269.4 & &  & 8,384  & 43,200 & \\
  $18A$ &  & 45 & 758 & &  6  & 49,315 & & 16,007 & 9,825.4 & &  & 16,007 & 43,200 & \\  
  $19A$ &  & 45 & 782 & &  6  & 41,403 & & 17,951 & 7,829.0 & &  & 17,951 & 43,200 & \\
  $20A$ &  & 45 & 674 & &  6  & 49,256 & & 8,676  & 9,495.3 & &  & 8,676  & 43,200 & \\
  \hline
\end{tabular}}
\end{center}
\end{table}

\subsection{MPSP-CRP with warm-start initialization}

The construction of the warm-start initialization includes the transcription of the heuristic solution into binary variables to input the model. The warm-started MPSP-CRP was solved using Gurobi\texttrademark{} 8.0 in a computer with an AMD Ryzen 7 1700 processor, 3 GHz, 64GB RAM, using 16 threads, and with Cplex\texttrademark{} 12.8 on a IBM cluster whose configuration consists of 8 nodes with 20 threads and 120GB RAM memory. The reason for this is to verify whether Cplex\texttrademark{} improves its performance using a computer with more threads and memory since it performed worse than Gurobi\texttrademark{} when using the same computer.

Tables \ref{tabResultadosMIPInstancias1} and \ref{tabResultadosMIPInstancias2} show the results of instances described in Tables \ref{tabInstancias1} and \ref{tabInstancias2}, respectively. The best result obtained by rule HR1 or HR2 for each instance was input as the initial solution for the solver. Following the pattern of the previous experiments, a limit of 432,000 seconds (5 days) was set for the computational time. After that, the solver is stopped and the result reached until then is reported. The first column, $I$, shows the number of the instance, the "Initial Solution" column indicates the objective function value of the solution obtained by the heuristics that is used as warm-start, the "Heuristic Time" column shows the time in seconds in which the initial solution was obtained by heuristics HR1 or HR2. Next, the $O.F.$ shows the value of the objective function obtained by the Cplex\texttrademark{} and Gurobi\texttrademark{} solvers, the $Time$ column shows the time in seconds that the respective solver took to solve the instance, and the $Gap$ column shows the gap between the solution found in relation to the best bound. The asterisk symbol ("$\ast$") means that the RAM memory available was not enough for the solvers to process these instances.

\begin{table}[h!]
\caption{Warm-start results: Cplex\texttrademark{} x Gurobi\texttrademark{}. \label{tabResultadosMIPInstancias1}}
\begin{center}
\scalebox{0.76}{
\begin{tabular}{|c|c|c|ccccc|ccccc|}
  \hline
  \multirow{2}{*}{$I$} &\multirow{2}{*}{\makecell{Initial \\Solution}} & \multirow{2}{*}{\makecell{Heuristic \\ Time}} & & \multicolumn{3}{c}{Cplex\texttrademark{}}  & & & \multicolumn{3}{c}{Gurobi\texttrademark{}}  &    \\
 \cline{5-7} \cline{10-12}
   & & & & $O.F.$ & $Time$ & $Gap$ & &  & $O.F.$ & $Time$ & $Gap$ &    \\
  \hline
  $1A$  & 1 & 0.13   & & 0  & 7.02       & 0.00\%  & & & 0 & 2.02       & 0.00\%  &  \\
  $2A$  & 4 & 0.14   & & 1  & 50.53      & 0.00\%  & & & 1 & 23.63      & 0.00\%  &  \\  
  $3A$  & 8 & 282.55 & & 2  & 13,052.05  & 0.00\%  & & & 2 & 5,927.86   & 0.00\%  &  \\
  $4A$  & 16& 804.34 & & 13 & 432,012.97 & 100\%   & & & 8 & 432,000.72 & 87.50\% &  \\
  \hline
  $5A$ & 0  & 0.11    & & 0  & 15.85       & 0.00\%   & & & 0  & 0.16       & 0.00\%  &  \\
  $6A$ & 5  & 0.16    & & 0  & 187,348.66  & 0.00\%   & & & 0  & 5,146.15   & 0.00\%  &  \\  
  $7A$ & 59 & 840.44  & & 59 & 432,000.72  & 89.53\%  & & & 59 & 432,001.78 & 83.05\% &  \\
  $8A$ & 48 & 1,339.2 & & 48 & 432,000.29  & 100.00\% & & & 48 & 432,001.02 & 98.41\% &  \\
  \hline
  $9A$  & 0  & 0.12   & &   0    & 1,544.49 & 0.00\% & & & 0      & 432.34 & 0.00\% &  \\
  $10A$ & 11 & 82.15  & & $\ast$ & $\ast$   & $\ast$ & & & $\ast$ & $\ast$ & $\ast$ &  \\  
  $11A$ & 26 & 218.74 & & $\ast$ & $\ast$   & $\ast$ & & & $\ast$ & $\ast$ & $\ast$ &  \\
  $12A$ & 68 & 562.40 & & $\ast$ & $\ast$   & $\ast$ & & & $\ast$ & $\ast$ & $\ast$ &  \\
  \hline
  \end{tabular}}
\end{center}
\end{table}

\begin{table}[h!]
\caption{Warm-start results: Cplex\texttrademark{} x Gurobi\texttrademark{}. \label{tabResultadosMIPInstancias2}}
\begin{center}
\scalebox{0.76}{
\begin{tabular}{|c|c|c|ccccc|ccccc|}
  \hline
  \multirow{2}{*}{$I$} &\multirow{2}{*}{\makecell{Initial \\Solution}} & \multirow{2}{*}{\makecell{Heuristic \\ Time}} & & \multicolumn{3}{c}{Cplex\texttrademark{}}  & & & \multicolumn{3}{c}{Gurobi\texttrademark{}}  &    \\
 \cline{5-7} \cline{10-12}
   & & & & $O.F.$ & $Time$ & $Gap$ & &  & $O.F.$ & $Time$ & $Gap$ &    \\
  \hline
  $1B$   & 0 & 0.11 & & 0 & 0.70     & 0.00\%  & & & 0 &  0.00   & 0,00\%       &  \\
  $2B$   & 0 & 0.12 & & 0 & 1.08     & 0.00\%  & & & 0 &  0.02   & 0,00\%       &  \\  
  $3B$   & 7 & 0.17 & & 0 & 194.55   & 0.00\%  & & & 0 & 255.71  & 0,00\%       &  \\
  $4B$   & 8 & 0.16 & & 0 & 4,228.83 & 0.00\%  & & & 0 & 2,548.9& 0,00\%       &  \\
  \hline
  $5B$  & 0 & 0.12    & & 0 & 15.91      & 0.00\%  & & & 0  & 0.17       & 0.00\%  &  \\
  $6B$  & 0 & 0.12    & & 0 & 25.27      & 0.00\%  & & & 0  & 0.25       & 0.00\%  &  \\  
  $7B$  & 40& 1,079.4 & & 40& 432,000.81 & 78.14\% & & & 11 & 59,820.72  & 0.00\%  &  \\
  $8B$  & 32& 0.35    & & 32& 432,001.41 & 96.87\% & & & 32 & 432,001.26 & 96.87\% &  \\
  \hline
  $9B$   & 0 & 0.16 & & 0      & 1,445.64 & 0.00\% & & &   0    & 674.13   & 0.00\% &  \\
  $10B$  & 0 & 0.26 & & 0      & 1,661.41 & 0.00\% & & &   0    & 1,003.85 & 0.00\% &  \\  
  $11B$  & 0 & 0.60 & & $\ast$ & $\ast$   & $\ast$ & & & $\ast$ & $\ast$ & $\ast$ &  \\
  $12B$  & 34& 1.48 & & $\ast$ & $\ast$   & $\ast$ & & & $\ast$ & $\ast$ & $\ast$ &  \\
  \hline
  \end{tabular}}
\end{center}
\end{table}

In general, by warm-starting the model, the solvers were able to find more solutions in less computational time. In instances $7A$, $8A$, $9A$, $8B$, $9B$ and $10B$, both solvers finished the optimization with a solution when using a warm-start. It is worth remembering that, without a warm-start initialization, those instances have reached the maximum computational time of 432,000 seconds and no feasible solution were found or processed at all.

It can be seen that when the warm-start is an optimal solution and there is enough memory, the optimization ends quickly, which was the case of instances $5A$, $9A$, $1B$, $2B$, $5B$, $6B$, $9B$, $10B$. Meanwhile, even though the computational time decreased in most instances, in the cases where the initial solution provided still has a large gap from the best bound, the solvers could take longer to reach the optimal solution (case of instances $2A$, $3A$ in Gurobi\texttrademark{} and $6A$ in Cplex\texttrademark{}), or not improve the initial solution at all (case of instances $7A$, $8A$, $8B$ and $7B$ in Cplex\texttrademark{}). Interestingly, in the case of instance $4A$, both solvers finished the optimization after 5 days with a solution worse than when a warm-start initialization was not used.

When comparing the performance between Cplex\texttrademark{} and Gurobi\texttrademark{} (Tables \ref{tabResultadosMIPInstancias1} and \ref{tabResultadosMIPInstancias2}), Gurobi\texttrademark{} finished the optimization in less or equal time in all instances than the Cplex\texttrademark{} (less in 13 instances), except for instances $6A$ and $3B$, and it reached a solution equal or better than Cplex\texttrademark{} in all instances (better in instances $4A$ and $7B$). Furthermore, in instances $4A$, $7A$, $8A$, even though both solvers could not improve the initial solution provided, Gurobi\texttrademark{} finished the optimization with a lower gap than the Cplex\texttrademark{}. 

In conclusion, the results showed that some limitations of the MPSP-CRP could be overcome within our warm-start approach, but only small to medium sized instances were exactly solvable, justifying the use of heuristic approaches to solve real-size instances. Additionally, even though Gurobi\texttrademark{} has performed better than the Cplex\texttrademark{}, it was notable that the difficulties in finding an optimal solution are not related with the solver selection but due to the problem features.







\section{Conclusions} \label{sec:conclusions}

This paper introduces an integrated optimization model for the multi-port stowage planning problem and the container relocation problem, called MPSP-CRP. The model was validated using computational experiments and a demonstration of the benefits of having an integrated model compared to solving each of the problems (CRP and MPSP) hierarchically was provided. Despite the commercial solvers being able to solve only small and medium sized instances, the model formulation represents an advance in the literature as it enables the comparison between the exact solution and the solution from the heuristics approaches, such as the HR1 and HR2, which provided optimal solutions for some small sized instances, and good solutions for medium, large and real sized instances within a short computational time.

Using heuristic solutions as a warm-start initialization in most of the cases has improved the speed of finding the solution by the exact method, but this improvement was not true for all cases. Therefore, it has provided a valuable understanding to be used in different real-world contexts. Analysing the results from all instances, one can also state that the yard relocation movement variables are most of the time equal to zero at the optimal solution, then it is possible to simplify the MPSP-CRP by setting this variable to zero in terms of searching for good quality solutions.

Extensions of the current research can be foreseen in three directions. One is related to the model representation itself, such as reformulating the model in a more simplified way, without losing its representation. Another direction would be to show more details from real contexts, such as including more ships, yards, different container sizes and weights and uncertainties regarding container arrivals and departures in ports. The last one would be to investigate several solution methods to solve large instances, such as those found in reality. Some methods could be evolutionary algorithms (e.g. Genetic Algorithms, Ant Colony, Bee algorithms, etc) and decomposition heuristics based on the model structure, or perhaps designing a Branch and Cut or Dynamic Programming framework to tackle the binary nature of the model more effectively. Finally, one day perhaps powerful computers such as those emerging in quantum computing could solve the model for large instances without no longer difficult.


\section*{Acknowledgments}

\noindent This research was supported by the S\~ao Paulo Research Foundation (FAPESP) (Process: 2015/24295-5), by the National Council for Scientific and Technological Development (CNPq) (Process: 400868/2016-4), and by IBM Company in Brazil. Moreover, this research benefited from a PhD grant from the Coordination for the Improvement of Higher Education Personnel (CAPES).

\section*{References}

\bibliography{bibliografiaESA}

\end{document}